\documentclass[12pt]{amsart}
\usepackage{amsmath}
 \usepackage{epsfig}
\usepackage{amssymb}
\usepackage{amsfonts}

\newtheorem{Proposition}{Proposition}
  \newtheorem{Remark}[Proposition]{Remark}
  
  \newtheorem{Lemma}[Proposition]{Lemma}
  
  \newtheorem{Theorem}[Proposition]{Theorem}

\newcommand{\be}[1]{\begin{equation}\label{#1}}
\newcommand{\ee}{\end{equation}}

\newcommand {\z}{{\noindent}}
\def\blackslug{\hbox{\hskip 1pt \vrule width 4pt height 8pt depth 1.5pt
\hskip 1pt}}
\def\qed{\quad\blackslug\lower 8.5pt\null\par}
 
\def\CC{\mathbb{C}}
 \def\RR{\mathbb{R}}
 
\def\ZZ{\mathbb{Z}}

\begin{document}

\ 

\vskip 2cm

\title{Failure of analytic hypoellipticity in a class of differential operators}


\gdef\shorttitle{Failure of analytic hypoellipticity in a class of PDOs}


%

\maketitle



{\large{ \begin{center}Ovidiu Costin, Rodica D.
      Costin\footnote{Research supported by NSF Grants 0074924, 0100495 and 0103807}\end{center}
}}

\

 \begin{center}
Department of Mathematics\\Busch Campus-Hill Center\\ Rutgers
       University\\ 110 Frelinghuysen Rd\\Piscataway, NJ 08854\\ 
e-mail: costin\symbol{64}math.rutgers.edu\\
\ \ \ \  rcostin\symbol{64}math.rutgers.edu 

\end{center}

\ 

\ 

\

\ 

{\bf{{ABSTRACT}}}  For the hypoelliptic differential operators $P={\partial^2_
    x}+\left( x^k\partial_ y -x^l{\partial_t}\right)^2$ introduced by
  T. Hoshiro, generalizing a class of M. Christ, in the cases of $k$
  and $l$ left open in the analysis, the operators $P$ also fail to be
  {\em{analytic}} hypoelliptic (except for $(k,l)=(0,1)$), in accordance
  with Treves' conjecture. The proof is constructive, suitable for
  generalization, and relies on evaluating a family of eigenvalues of
  a non-self-adjoint operator.



\vskip 1cm

\section{Introduction and result}

A differential operator $P$ is said to be hypoelliptic (respectively,
analytic hypoelliptic) on $\Omega$ if for any $C^\infty$
(respectively, $C^\omega$) function $f$ on some open set
$U\subset\Omega$ all the solutions $u$ of $Pu=f$ belong to
$C^{\infty}(U)$ (respectively, $C^{\omega}(U)$).

The basic result about the hypoellipticity of operators of the type
``sum of squares'', $P=X_1^2+...+X_n^2$ , where $X_1,...,X_n$ are real
vector fields of class $C^\omega(\Omega)$, is H\"ormander's theorem
\cite{Hormander} which gives necessary and sufficient conditions for
hypoellipticity. But these assumptions are not sufficient for analytic
hypoellipticity as was first proved by Baouendi and Goulaouic
\cite{BG}. Other classes of hypoelliptic operators which fail to be
analytic hypoelliptic have been found and there are important results
on analytic regularity (see Christ \cite{Christ2}-\cite{C1996}, Christ
and Geller \cite{C-G}, Derridj and Tartakoff
\cite{D-T88}-\cite{D-T95}, Derridj and Zuily \cite{D-Z}, Francsics and
Hanges \cite{F-H88}-\cite{F-HDuke}, Grigis and Sj\" ostrand
\cite{G-S}, Hanges and Himonas \cite{H-H91}-\cite{H-H98}, Helffer
\cite{Hlf}, Hoshiro \cite{Hoshiro}, Metivier \cite{M80}-\cite{M81},
Pham The Lai and Robert \cite{L-R}, Sj\" ostrand \cite{S83} ,Tartakoff
\cite{T80} ,Treves \cite{T78}-\cite{Treves}, see also the survey
\cite{F-HSurvey} for more references).

However, the question of finding a general characterisation of
analytic hypoellipticity for sum of squares operators is still open.
Christ gave a criterion for analytic hypoellipticity in the
two\--dimensional case \cite{C1996}. Treves conjectured a general
criterion for analytic hypoellipticity \cite{Treves}: ``For a sum of
squares of analytic vector fields to be analytic hypoelliptic it is
necessary and sufficient that every Poisson stratum (defined by the
symbols of the vector fields) of its characteristic variety be
symplectic.''

Extending a result of Christ \cite{Christ2}, Hoshiro gave new
examples of hypoelliptic operators in $\RR^3$ which fail to be
analytic hypoelliptic \cite{Hoshiro}. They have the form

\begin{equation}\label{opehoshi}
P=\frac{\partial^2}{\partial x^2}+\left( x^k\frac{\partial}{\partial
    y} -x^l\frac{\partial}{\partial t}\right)^2
\end{equation}
where $k<l$ are non-negative integers. It is shown, through an elegant
proof, that $P$ is not analytic hypoelliptic if either of the following
assumptions are satisfied:

(i) $(l+1)/(l-k)$ is not a positive integer.

(ii) Both $l-k$ and $(l+1)/(l-k)$ are odd integers.

The case $k=0$ was studied by Christ \cite{Chr01} and further refined by
Yu \cite{Yu}.  In view of Treves' conjecture, it is interesting to
investigate the cases remained open.

The purpose of the present paper is to show that in the remaining
cases, except for $(k,l)=(0,1)$, the operator $P$ fails to be analytic
hypoelliptic as well. 

This result is in agreement with Treves' conjecture.  Indeed, the
Poisson strata of the characteristic variety
$\Sigma=\{\xi=x^k\eta-x^l\tau=0\}$ of the operators (\ref{opehoshi})
with $k<l$ are as follows. For $k\geq 1$ then $\Sigma_0=\{
\xi=x^k\eta-x^l\tau=0, x\ne 0\}$, $\Sigma_1=\{ x=\xi=0,\eta\ne 0\}$,
and $\Sigma_2=\{ x=\xi=\eta=0,\tau\ne 0\}$ which is not symplectic
since its codimension is odd. For $k=0$ and $l\geq 2$ the strata are
$\Sigma_0=\{ \xi=\eta-x^l\tau=0, x\ne 0\}$ and $\Sigma_1=\{
x=\xi=\eta= 0,\tau\ne 0\}$ which is not symplectic due to its odd
codimension. By contrast, for $k=0$ and $l=1$ the characteristic
variety is symplectic, and the operator is analytic hypoelliptic (see
also \cite{Christ2}).

\section{Outline of proof}

The proof uses a standard reduction to ordinary differential equations:

\begin{Lemma}[\cite{Hoshiro}]\label{lemahoshi}

Let, for $\zeta\in\CC$, 

\begin{equation}\label{Pzeta}
P_\zeta=-\frac{d^2}{dx^2}+(x^l-x^k\zeta)^2
\end{equation}

\z If there exist $\zeta\in\CC$ and $f\in L^\infty(\RR)$, $f\ne
0$, satisfying 

\begin{equation}\label{eqfx}
P_\zeta f=0
\end{equation}

\z then the operator $P$ is not analytic hypoelliptic.

\end{Lemma}

The proof of Lemma \ref{lemahoshi} and further references are found in
\cite{Hoshiro}.

The results of \cite{Hoshiro} and of the present paper give:

\begin{Theorem}\label{failureofah}
  
  For all non-negative $k,l$ with $k<l$ and $(k,l)\ne (0,1)$ there are
  values of $\zeta\in\CC$ such that equation (\ref{eqfx}) has a
  solution $f\in L^{\infty}(\RR)$, $f\not\equiv 0$.

As a consequence (\ref{opehoshi}) is not analytic hypoelliptic.

\end{Theorem}

To complete the proof of Theorem \ref{failureofah} we only need to
consider $k$ and $l$ which fail both conditions (i) and (ii), and we prove:

\begin{Proposition}\label{fmic}
  
  For all integer values of $k,l$ with $1\leq k<l$ which fail both
  conditions (i) and (ii) there are values of $\zeta\in\CC$ such that
  equation (\ref{eqfx}) has a solution $f\in L^{\infty}(\RR)$,
  $f\not\equiv 0$.

An array of such values of $\zeta$ satisfy the estimate

\begin{equation}\label{valzeta}
\zeta=\zeta_M=e^{\frac{i\pi(l-k)}{2(l+1)}}\left[
  M\frac{\pi(l+1)(k+1)}{l-k}\right]^{\frac{l-k}{l+1}}\, \left[
  1+O\left(M^{-1}\right)\right]\ \ \ (M\rightarrow +\infty)
\end{equation}

\z for $M\in\ZZ_+$ large enough.

\end{Proposition}
 
Note that $k=0$ is not considered in Proposition \ref{fmic}. In fact
for $k=0$ and $l\geq 2$ condition (i) is satisfied so failure of
analytic hypoellipticity was proved in \cite{Hoshiro}.  The case $k=0$
was studied earlier in \cite{Christ2} (see also the references
therein) where it was shown that such operators are not analytic
hypoelliptic except for the case $l=1$. The case $(k,l)=(0,1)$ is
fundamentally distinct as in this case the operator (\ref{opehoshi})
is analytic hypoelliptic \cite{Christ2}; this distinction appears
naturally in our proof (see Remark \ref{k0l1}).

Our proof of existence of bounded solutions of (\ref{fmic}) differs
from the approaches in \cite{Christ2} and \cite{Hoshiro} and the
method is constructive and general; we restrict the proof to the cases
not covered in \cite{Hoshiro} for simplicity. The values of $\zeta$ of
Theorem \ref{failureofah} are related to the eigenvalues of a
non-selfadjoint operator associated to (\ref{Pzeta}). Allowing $\zeta$
to be large is a crucial element in our proof, as it leads to
substantial asymptotic simplification of the equations involved.
There are other arrays of eigenvalues $\zeta$ besides (\ref{valzeta}),
but we will not pursue this issue here.

\section{Main Proofs}

\subsection{Notations}\label{Notations}

We denote by $const\, $ positive constants independent of the
parameter $\xi$.

$O(\rho^{-c})$ will denote functions (of the variable $v$, or of $s$)
, depending on the parameter $\rho$, which are less, in absolute
value, than $const\, \rho^{-c}$ uniformly on specified paths.

$o(1)$ will stand here only for terms decaying like powers of $\rho$,
i.e. for $O(\rho^{-c})$ for some $c>0$.

\subsection{General setting}

Classical methods in the theory of ordinary differential equations
were used to show:

\begin{Lemma}[\cite{Hoshiro}]\label{L1Hoshi}
  Equation (\ref{eqfx}) has solutions which decrease exponentially
for $x\rightarrow +\infty$, respectively $x\rightarrow -\infty$ (see \cite{Hoshiro}):

\begin{equation}\label{decrfx}
\begin{array}{ll}
f_-^{\{+\infty\}}(x)\sim x^{-\frac{l}{2}}\exp\left({-W(x)}\right)\ \ \ (x\rightarrow
+\infty)\\
f_-^{\{-\infty\}}(x)\sim x^{-\frac{l}{2}}\exp\left({(-1)^lW(x)}\right)\ \ \ (x\rightarrow
-\infty)
\end{array}
\end{equation}

\z where $W(x)=x^{l+1}/(l+1)-\zeta x^{k+1}/(k+1)$, and solutions with
exponential growth

\begin{equation}\label{incrfx}
\begin{array}{ll}
f^{\{+\infty\}}_+(x)\sim x^{-\frac{l}{2}}\exp\left({W(x)}\right)\ \ \ (x\rightarrow
+\infty)\\
f^{\{-\infty\}}_{+}(x)\sim x^{-\frac{l}{2}}\exp\left({(-1)^{l+1}W(x)}\right)\ \ \ (x\rightarrow
-\infty)
\end{array}
\end{equation}

\end{Lemma}

We will show that there are values of $\zeta$ for which
$f^{\{+\infty\}}_{-}$ is a multiple of $f^{\{-\infty\}}_-$ by finding the 
analytic continuation of $f^{\{+\infty\}}_-$ from $+\infty$ to $-\infty$
along $\RR$.

{\em{Notation.}} To simplify the notations we will omit from now on
the superscript $+\infty$ from the notations of the small, or a large
solution for $x\rightarrow +\infty$; we however need to keep the
superscript $-\infty$ to make the distinction for the small, or a
large solution as $x\rightarrow -\infty$.

Denote $m=l-k\ge 1$ and $q=(l+1)/(l-k)-2$. Since we assume that
$(l+1)/(l-k)$ is a positive integer, then $q\geq 0$, $q\in\ZZ$.
If $m$ is odd we assume that $q$ is even (since the case with
$m$ and  $q$ odd falls under the already studied case (ii) \cite{Hoshiro}).

The substitution 

\begin{equation}\label{vzx}
v=\zeta^{-1}x^{l-k}
\end{equation}

\z $g(v)=f(x)$, brings (\ref{eqfx}) to the form

\begin{equation}\label{eqgv}
g''(v)+\frac{m-1}{m}\frac{1}{v}g'(v)=-\xi^2v^{2q}\left( v-1\right)^2g(v)
\end{equation}

\z where
 
\begin{equation}\label{defxi}
\xi=-\frac{i}{m} \zeta^{q+2}
\end{equation}

We will only consider $\xi$ large enough in a strip containing $\RR_+$
(cf. (\ref{valzeta})):

\begin{multline}\label{orderz}
  \xi=\rho+iZ\ \ {\mbox{with}}\ \rho\in\RR_+\ \ {\mbox{large\ enough,}}\\
  {\mbox{and\ }}Z\in\RR\, ,\ \  Z=O(1)\ \  (\rho\rightarrow +\infty)\ 
\end{multline}

\z(see \S\ref{Notations} for notations).

Correspondingly,
 
\begin{equation}\label{argzeta}
\zeta=e^{i\frac{\pi}{2(q+2)}} (\rho m)^{\frac{1}{q+2}} \left[ 1+o(1)\right]\ \ \
\ (\rho\rightarrow +\infty)
\end{equation}

The path $-\RR$ from $+\infty$ to $-\infty$ in the $x$-plane
corresponds to the path $v(-\RR)$ on the Riemann surface above
$\CC\setminus\{0\}$ of the $v$-plane.  If $m$ is even then $v(-\RR)$
comes from $\infty$ along the half-line
$$d = e^{-i\frac{\pi}{2(q+2)}} \RR_+$$ then turns by an angle of $m\pi$
around $v=0$ and  returns towards $\infty$ along $d$. If $m$ is odd then
$v(-\RR)$ comes from $\infty$ along $d$, turns by an angle $m\pi$ around
$v=0$ then follows $-d$ towards $\infty$.

To prove Proposition \ref{fmic} we show that there are values of $\xi$
satisfying (\ref{orderz}) for which equation (\ref{eqgv}) has a
bounded solution on $v(-\RR)$ (in fact, exponentially small for
$v\rightarrow\infty$ on the path).

The proof is done as follows. For not too small $v\in d$ we find
precise estimates of the unique (up to a multiplicative factor) small
solution $g_-(v)$ of (\ref{eqgv}) (Lemma \ref{L1}). On the other hand,
for small $v$, we find precise estimates for two independent solutions of
(\ref{eqgv}), $g_a(v)$ and $g_s(v)=v^{1/m}g_{a;s}(v)$ with
$g_a,g_{a;s}$ analytic (Lemma \ref{L2}). The two regions where the said
estimates hold overlap (for $v$ small, but not too close to $0$),
so the two representations can be matched in this common region
((\ref{firstmatch}), Lemma \ref{asygas}).

If $m$ is even, we find values of $\xi$ for which $g_-$ is a multiple
of $g_a$; then the continuation of $g_-$ along $v(-\RR)$ yields a
small solution (Lemma \ref{Csis0}) and Proposition \ref{fmic} is
proved for $m$ even.

If $m$ is odd, and $q$ is even, we also find precise estimates of the
small solution $g_{-}^{\{-\infty\}}(v)$ of (\ref{eqgv}) for $v\in -d$, from
$\infty$ down to not too small $v$ (Lemma \ref{L1min}). We decompose $g_-$
along $g_a$, $g_s$ ((\ref{firstmatch}), Lemma \ref{asygas}), perform
analytic continuation upon a rotation by angle of $m\pi$ then determine
$\xi$ so that this result is a multiple of $g_{-}^{\{-\infty\}}(v)$ 
((\ref{decmin}), Lemma \ref{L10p9}). This ends the proof of Proposition
\ref{fmic}.

\subsection{Small solutions of (\ref{eqgv}) for $|v|\rightarrow
  \infty$, $v\in d$}\label{smallsol}

In \S\ref{smallsol} we find precise estimates for $g_-$, the exponentially
small solution of (\ref{eqgv}) on $d$; the main result
is Lemma \ref{L1} (ii) and (\ref{ttxxvp}).

\subsubsection{Setting}\label{settinggm}
We make suitable substitutions in (\ref{eqgv}) to obtain uniform
estimates for a relatively large domain of variation of $v$.

Let $a$ be a number satisfying

\begin{equation}\label{vala}
\frac{1}{q+1}\left[ 1-\frac{2m}{3m(q+1)+2m-1}\right]<a<\frac{1}{q+1}
\end{equation}

Note that $m(q+1)=k+1>1$ under our assumptions on $k,l$. Also, a
simple calculation shows that (\ref{vala}) implies $a>1/(q+2)$.

Denote $s(v)=i\xi\int_1^v t^q(t-1)dt$; then

\begin{equation}\label{sdev}
s(v)=i\xi\left(\frac{1}{q+2}v^{q+2}-\frac{1}{q+1}v^{q+1}+P_0\right) 
\equiv i\xi(v-1)^2P(v)
\end{equation}

\z We note that $P(v)=(q+1)^{-1}(q+2)^{-1}\sum_{j=0}^q(j+1)v^j$ hence
$P(1)\ne 0$ and
$$P_0= P(0)=(q+1)^{-1}(q+2)^{-1}$$

The inverse function $v(s)$ of (\ref{sdev}) is algebraic, ramified only
when $v=1$ (at $s=0$) and $v=0$ (at $s=s_0\equiv i\xi P_0$, if $q\ne
0$).

The substitution 

\begin{equation}\label{subgg1}
g(v)=P(v)^{1/4}v^{q_1}g_1(v)\ \ \ \ \ \ \ \ \ 
{\mbox{where\ \ \ \ }}q_1=-\frac{q+1}{2}+\frac{1}{2m}
\end{equation}

\z followed by

\begin{equation}\label{number}
s=s(v)\ \ ,\ \ g_1(v)=h\left(s(v)\right)
\end{equation}

\z transform (\ref{eqgv}) into

\begin{equation}\label{eqhs}
h''(s)+\frac{1}{2s}h'(s)-h(s)=R(s)h(s)
\end{equation}

\z where

\begin{equation}\label{Rdef}
R\left(s\right)=\frac{\tilde{R}(v)}{\xi^2v^{2q+2}(v-1)^2}\, \Big| _{v=v(s)}
\end{equation}

\z where

\begin{equation}\label{defRtilde}
\tilde{R}(v)=-\frac{3}{16}\left(\frac{vP'(v)}{P(v)}\right)^2+\frac{1}{4}\frac{v^2P''(v)}{P(v)}-\frac{q}{4}\frac{vP'(v)}{P(v)}+\left(\frac{q+1}{2}\right)^2-\frac{1}{4m^2}
\end{equation}

The branch of $v(s)$ in (\ref{Rdef}) is $v\in d$ if $s\in s(d)$,
respectively $v\in -d$ if $s\in s(-d)$; we denote the latter by
$v^-(s)$.

\subsubsection{Shifting the path $s(d)$}\label{2.1.1.}

We will use some known facts on the behavior of a solution of a
differential equation towards a rank one irregular singular point to
calculate solutions of (\ref{eqgv}) on a path different from $s(d)$.

By Lemma~\ref{L1}, equation (\ref{eqgv}) has an exponentially small
solution $g_-$ along $d$, unique up to a multiplicative factor (since
$g_-(v)=f_-^{\{+\infty\} }(x)$, see (\ref{decrfx}), (\ref{vzx})).

After the substitution (\ref{subgg1}), (\ref{number}) equation
(\ref{eqgv}) becomes (\ref{eqhs}). The coefficients of the linear
equation (\ref{eqhs}) are analytic on the Riemann surface above
$\CC\setminus\{0, s_0\}$.  In addition, $s=0$ is an algebraic
singularity for the solutions of (\ref{eqhs}). For large values of $s$
on one sheet of this surface standard theory of linear equations
applies.\footnote{The study of solutions of differential equations
  (linear or nonlinear) having a prescribed asymptotic behavior is done
  on shifted sectors, for large values of the independent variables
  (for linear equations see \cite{Sibuya}, \cite{Varadarajan}, and for
  nonlinear equations see \cite{Wasow}).}

We need the small solution $h_-$ of (\ref{eqhs}) on $s(d)$ as we can then
define $g_-(v)=P(v)^{1/4}v^{q_1}h_-(s(v))$ to be the small solution of
(\ref{eqgv}) on $d$.

The path $s(d)$ comes from $\infty$ in the right half-plane, crosses
the left half-plane and ends at $s_0$--- see the Appendix \S\ref{s(d)}
for details.  By standard theory of linear equations (see
\cite{Sibuya}, \cite{Wasow}) there is a solution $h_-$ of (\ref{eqhs})
which is exponentially small for $|s|\rightarrow +\infty$, $\arg(s)\in
(-\frac{\pi}{2},\frac{\pi}{2})$, and its asymptotic representation has
the form

\begin{equation}\label{asyhmin}
h_-(s)=s^{-1/4}e^{-s}(1+O(s^{-1}))
\end{equation}

\z which is valid for $\arg(s)\in (-\frac{\pi}{2},\pi)$ ($|s|\rightarrow
\infty$). Similarly, there exists a solution satisfying

\begin{equation}\label{asyhplus}
h_+(s)=s^{-1/4}e^{s}(1+O(s^{-1}))
\end{equation}

\z for $|s|\rightarrow +\infty$, $\arg(s)\in (-\pi/2,\pi)$.\footnote{The
  intervals stated for the validity of asymptotic representations are
  not optimal, but we do not need more.}

We will follow the values of $g_-(v)$ along $d$ from $\infty$ down to
small values of $v$: 

\begin{equation}\label{vmmnou}
v_{P}=e^{-i\frac{\pi}{2(q+2)}} \rho^{-a}
\end{equation}

\z on 

\begin{equation}\label{da}
d_a\equiv v_{P} [1,+\infty)\subset d
\end{equation}

In variable $s$ this amounts to following $h_-(s)$ along $s(d)$ from
$\infty$ down to

\begin{equation}
s_{P}=s\left( v_{P}\right)=i\rho
P_0-e^{i\frac{\pi}{2(q+2)}}\frac{\rho^{1-a(q+1)}}{q+1}-ZP_0+o(1)\ \
\ \ (\rho\rightarrow +\infty)
\end{equation}

The behavior (\ref{asyhmin}) of $h_-$ (respectively, (\ref{asyhplus})
of $h_+$) holds for $s\rightarrow \infty$ , $s\in s(d)$ and also for
$s\rightarrow \infty$ , $s\in i\RR_+$. It is more convenient to use
the latter path to continue $h_\pm$ from $\infty$ to the finite plane.

We define a path $\ell$ coming from infinity along $i\RR_+$ and ending
at $s_{P}$, homotopically equivalent to $s(d_a)$ in a subregion of
$\CC\setminus \{0,s_0\}$ in which $\arg
s\in(-\frac{\pi}{2},\frac{\pi}{2}]$ for $s$ large. Then the value of
$h_-$ at $s_{P}$ can be calculated along $\ell$.

Let $t_0>0$ small enough, independent of $\rho$. Let 

\begin{equation}\label{sm}
s_M=i\rho P_0-i\frac{\rho^{1-a(q+1)}}{q+1}+it\ \ \ ,\
t\in[-t_0,0]\subset\RR_-
\end{equation}

\z with $s_M=s(v_M)$ where 

\begin{multline}\label{vm}
v_M=\rho^{-a}\left( 1+(iZP_0-t)\rho^{-1+a(q+1)}(1+\Delta)\right)\ ,\\
{\mbox{with}}\ \  \Delta=O(\rho^{-a})\ (\rho\rightarrow +\infty)
\end{multline}

\begin{figure}
\epsfig{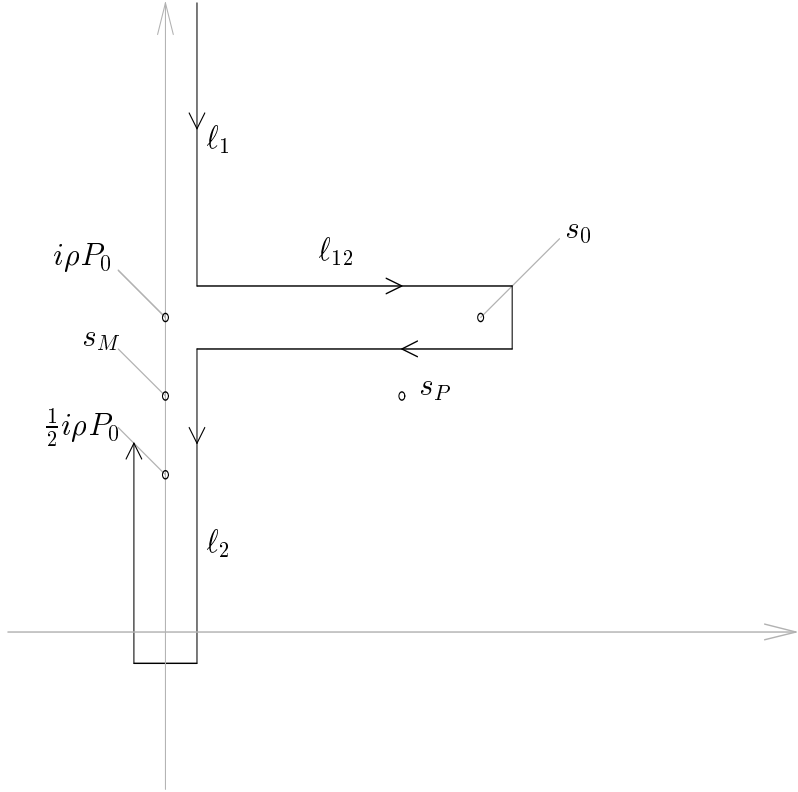}
\\
\caption{The beginning of the path $\ell$: $\ell_1$, $\ell_{12}$ and $\ell_2$.}
\end{figure}

Let $\ell $ be the path in the $s$-plane consisting of the following
five parts. (Figure 1 sketches the first three parts of $\ell$ in the
case $\Re s_0\geq 0$; the vertical lines are in fact included in
$i\RR_+$.) Let $\ell_1$ be the half-line coming from infinity along
$i\RR_+$ down to $i\rho P_0$. If $\Re s_0\geq 0$ let also $\ell_{12}$
be the segment from $i\rho P_0$ to $s_0$, avoiding $s_0$ clockwise on
a small circle, then returning to $i\rho P_0$. (If $\Re s_0< 0$ then
$\ell_1$ continues with $\ell_2$.) Let $\ell_2$ come from $i\rho P_0$
along $i\RR_+$ towards $0$, avoiding $0$ clockwise on a small circle,
then go up to $i\rho P_0/2$. Let $\ell_3$ be the part from $i\rho
P_0/2$ to $s_M$. Finally, $\ell_4$ is the segment $[s_M,s_{P}]$.

The branch of $v(s)$ in (\ref{Rdef}) on $\ell$ is obtained by analytic
continuation of the branch on $s(d_a)$ along the path that gives the
homotopy between $s(d_a)$ and $\ell$ (see \S\ref{Someestim} for details).

\subsubsection{The solution $g_-$}

The solutions of the left-hand side of (\ref{eqhs}), i.e. which satisfy 
\begin{equation}\label{eqB}
B''(s)+\frac{1}{2s}B'(s)-B(s)=0
\end{equation}
will turn out to be very close to solutions of (\ref{eqhs}) on $\ell$.

The solutions of (\ref{eqB}) are Bessel functions, with two
independent solutions given by

\begin{equation}\label{formBpm}
B_\pm(s)=\int_1^{+\infty}e^{\pm ps}\left( p^2-1\right)^{-3/4}\, dp\ \
\ \ \ \mbox{for}\ \  s\in i\RR_+
\end{equation}

Equation (\ref{eqhs}) can be written in integral form as

$$h(s)=C_1B_+(s)+C_2B_-(s)+\mathcal{J}(h)(s) $$

\z where

\begin{equation}\label{defJ}
\mathcal{J}(h)(s)={W_0}^{-1}\int_{+i\infty}^s\left[ B_+(s)B_-(\sigma)-
B_-(s)B_+(\sigma)\right]{\sigma^{1/2}}R(\sigma)h(\sigma)\,
d\sigma
\end{equation}

\z $W_0$ is a constant independent of $\xi$ (the Wronskian of $B_+,B_-$)
and the path of integration is along $\ell$.

For $C_1=0$ and $C_2=1$ we get the small solution (while for $C_2=0$ and
$C_1=1$ we get an independent, large solution):

\begin{Lemma}\label{L1}

(i) The integral equation

\begin{equation}\label{eqhmin}
h_-(s)=B_-(s)+\mathcal{J}(h_-)(s)
\end{equation}

\z has a unique solution on $\ell$ and 

\begin{equation}\label{ttxxvp}
\sup_{s\in\ell}|h_-(s)-B_-(s) |<const\, 
\rho^{-\frac{5}{4}+a(q+1)}
\end{equation}

(ii) The solution $g_-(v)=P(v)^{1/4}v^{q_1}h_-(s(v))$ of
(\ref{eqgv}) tends to zero for $v\rightarrow \infty$, $v\in d$.

(iii) $g_-$ depends analytically on $\xi$ for $\xi$
large in a strip containing $\RR_+$.

\end{Lemma}

The proof is given in \S\ref{ProofL1}.

\subsection{Solutions of (\ref{eqgv}) for $|v|\rightarrow
  \infty$, $v\in -d$ when $m$ is odd and $q$ is even}\label{smallsolm}

In this section we find precise estimates for $g_-^{\{-\infty\}
  }$, the exponentially small solution of (\ref{eqgv}) on $ -d$, and
for the large solution $g_+^{\{-\infty\} }$ and the main result is
Lemma \ref{L1min} (ii) and (\ref{ttxxvm}). The construction is similar
  to the one for $v\in d$; we provide the details for completeness. 

The proof uses constructions similar to those of \S\ref{smallsol}.

\subsubsection{Shifting the path $s(-d)$}\label{shiftsmd}

For $q$ even $s(-d)$ starts at $s_0$, then goes to $\infty$ in the
right-half plane (see Appendix \S\ref{s(d)} for details).

We consider $g_-^{\{-\infty\} }$ along $-d_a$ (see (\ref{da}),
(\ref{vmmnou})). Denote

\begin{equation}\label{sPm}
s_{P}^-=s\left(- v_{P}\right)=i\rho
P_0+e^{i\frac{\pi}{2(q+2)}}\frac{\rho^{1-a(q+1)}}{q+1}-ZP_0+o(1)\ \
\ \ (\rho\rightarrow +\infty)
\end{equation}

The behavior (\ref{asyhmin}), respectively (\ref{asyhplus}) of the
solution $h_-^{\{-\infty\} }$, respectively $h_+^{\{-\infty\} }$, holds for
$s\rightarrow \infty$ , $s\in s(-d)$ and also for $s\rightarrow
\infty$ , $s\in i\RR_+$. It is more convenient to use the latter path
to continue $h^{\{-\infty\}}_{\pm}$ in the finite plane.

We define a path $\ell^-$ coming from $\infty$ along $i\RR_+$ and
ending at $s_{P}^-$, homotopically equivalent to $s(-d_a)$ in a subset of
$\CC\setminus \{0,s_0\}$ for which $\arg s\in
(-\frac{\pi}{2},\frac{\pi}{2}]$ 
for $s$ large. Then the value of
$h^{\{-\infty\}}_{\pm}$ at $s_{P}^-$ can be calculated along $\ell^-$.

Let

\begin{equation}\label{smmin}
s_M^-=i\rho P_0+i\frac{\rho^{1-a(q+1)}}{q+1}+it\ \ \ ,\ t\in[0,t_0]\subset\RR_+
\end{equation}

\z (with $t_0$ small). We have $s_M^-=s(v_M^-)$ for

\begin{multline}\label{vmmin}
v_M^-=v^-\left( s_M^-\right) \\
=-\rho^{-a}\left( 1-(iZP_0-t)\rho^{-1+a(q+1)}(1+\Delta_-)\right)\ , \ \ \ \Delta_-=O(\rho^{-a})\ (\rho\rightarrow +\infty)
\end{multline}

Let $\ell^-$ be the path in the $s$-plane consisting of $(+i\infty,,
s_M^-]\subset i\RR_+$ followed by the segment $[s_M^-,s_{P}^-]$.

The branch of $v(s)$ in (\ref{Rdef}) on $\ell^-$ is obtained by analytic
continuation of the branch $v^-(s)$ on $s(-d_a)$ along the path that gives the
homotopy between $s(-d_a)$ and $\ell^-$.

\subsubsection{The solutions $g_{\pm}^{\{-\infty\}}$}

\begin{Lemma}\label{L1min}

(i) The integral equations

\begin{equation}\label{eqhpmvm}
h_{\pm}^{\{-\infty\} }(s)=B_\pm(s)+\mathcal{J}(h_{\pm}^{\{-\infty\} }(s)
\end{equation}

\z have unique solutions on $\ell^-$ and 

\begin{equation}\label{ttxxvm}
\sup_{s\in\ell_-}|h_{\pm}^{\{-\infty\} }(s)-B_\pm(s) |<const\, \max\{
\rho^{-3/4},\rho^{-1+a(q+1)}\}
\end{equation}

(ii) The solution
$g_{-}^{\{-\infty\} }(v)=P(v)^{1/4}v^{q_1}h_{-}^{\{-\infty\} }(s(v))$ of
(\ref{eqgv}) tends to zero for $v\rightarrow \infty$, $v\in - d$ and
$g_{+}^{\{-\infty\} }(v)=P(v)^{1/4}v^{q_1}h_{+}^{\{-\infty\} }(s(v))$ tends to
infinity for $v\rightarrow \infty$, $v\in - d$.

(iii) The solutions $g_{\pm}^{\{-\infty\} }$ depend analytically on
$\xi$ for $\xi$ large, in a strip containing $\RR_+$.

\end{Lemma}

The proof is given in \S\ref{ProofL1m}.

\subsection{Solutions for small $v$}\label{anbrsol}

The linear equation (\ref{eqgv}) has a regular singular point at
$v=0$; it has an analytic solution $g_a$ and a ramified solution
$g_s=v^{1/m}g_{a;s}$ (with $g_{a;s}$ analytic). For a proof and estimates we
separate the dominant terms (for small $v$), writing (\ref{eqgv}) in
the form

\begin{equation}\label{eqgvsmall}
g''(v)+\frac{m-1}{m}\frac{1}{v}g'(v)+\xi^2v^{2q}g(v)=-\xi^2v^{2q}\left(v^2-2v\right)g(v)
\end{equation}

The the LHS (dominant part) of (\ref{eqgvsmall}) has two independent
solutions (see \S\ref{Wnudepxi} for details)

\begin{multline}\label{gags}
\tilde{g}_a(v)=\int_{-1}^1 e^{-ip\xi\frac{v^{q+1}}{q+1}}\left(
  1-p^2\right)^A\, dp\ \ ,\ {\mbox{with}}\ A=-\frac{1}{2}\left(
  1+\frac{1}{m(q+1)}\right)\\
\tilde{g}_s(v)=v^{1/m}\int_{-1}^1 e^{-ip\xi\frac{v^{q+1}}{q+1}}\left(
  1-p^2\right)^B\, dp\ ,\ {\mbox{with}}\ B=-\frac{1}{2}\left( 1-\frac{1}{m(q+1)}\right)
\end{multline}

\z We have $A,B\in (-1,0)$ since $m(q+1)> 1$.

\begin{Remark}\label{k0l1}
  
  The condition $A,B\in (-1,0)$ is fundamental to our proof of
  Proposition \ref{fmic}, as the behavior of solutions $g(v)$ for
  small $v$ is different in other cases (cf. Remark \ref{aint1inf}).
  The only pair of values of the nonnegative integers $k<l$ for which this
  condition fails is $(k,l)=(0,1)$. It is known \cite{C1996} that in
  this case $P$ is in fact analytic hypoelliptic.

\end{Remark}

Then equation (\ref{eqgvsmall}) can be written in integral form (see
\S\ref{Wnudepxi} for details)

\begin{equation}\label{ginteq}
g(v)=C_1\tilde{g}_a(v)+C_2\tilde{g}_s(v)+\mathcal{G}(g)(v)
\end{equation}

\z where 

\begin{equation}\label{defcalG}
\mathcal{G}(g)(v)=const\, \xi^{r+1}\int_0^v
\left[ \tilde{g}_a(v)\tilde{g}_s(t)-\tilde{g}_s(v)\tilde{g}_a(t)\right]t^{2(q+1)-\frac{1}{m}}(t-2)g(t)\,dt 
\end{equation}

\z where 

\begin{equation}\label{defr}
r=1-\frac{1}{m(q+1)}
\end{equation}

Let $\mathcal{S}$ be the path in the $v$-plane composed of two
segments: $\mathcal{S}_1=[0,v_F]$ with $v_F=\kappa_0\rho^{-a}$ (where
$\kappa_0^{q+1}=(1+z)^{-1}$, $\Re\kappa_0>0$) followed by the segment
$\mathcal{S}_2=[v_F,v_{M}]$.

\begin{Lemma}\label{L2}
  
  The integral equation (\ref{ginteq}) has a unique solution for $v\in
  \mathcal{S}$.  For $C_2=0$, $C_1=1$, its solution $g_a(v)$ is
  analytic at $v=0$.  For $C_1=0$, $C_2=1$, its solution has the form
  $g_s(v)=v^{1/m}g_{a;s}(v)$ with $g_{a;s}(v)$ analytic.

The following estimates hold

\begin{equation}\label{estimgasD}
\left| g_{a,s}(v)-\tilde{g}_{a,s}(v)\right| < const\,
\rho^{1-a(q+2)}\ \ {\mbox{for}}\ v\in\ \mathcal{S} 
\end{equation}

Also, $g_a$, $g_s$ are analytic in $\xi$ for $\xi$ large in a strip containing $\RR_+$.

\end{Lemma}

The proof in given in \S\ref{ProofL2}.

\subsection{Matching between $g_-$ and $g_a,g_s$}\label{locmat}

The calculations are done at $v=v_M$, see (\ref{vm}). Note that
$v_M$ depends on the real parameter $t$ (independent of $\xi$).

We first need asymptotic estimates for the solutions $g_a,g_s,g_-$:

\begin{Lemma}\label{asygminus}

For $v$ satisfying (\ref{vm}) we have for $\xi\rightarrow\infty$

\begin{equation}\label{aga}
g_{a}(v)=K_{a}\rho^{-\beta}\left[e^{i(R+\lambda)}+e^{i\pi(1+A)}e^{-i(R+\lambda)}\right]+O(\rho^{-d})
\end{equation}

\begin{equation}\label{ags}
g_{s}(v)=K_{s}\rho^{-\beta}\left[e^{i(R+\lambda)}+e^{i\pi(1+B)}e^{-i(R+\lambda)}\right]+O(\rho^{-d})
\end{equation}

\z and

\begin{equation}\label{agm}
g_-(v)=K\rho^\gamma\left[
e^{iP_0\xi}e^{-i(R+\lambda)}+\frac{i}{\sqrt{2}}e^{-iP_0\xi}e^{i(R+\lambda)}+O(\rho^{-1+a(q+1)})\right]
\end{equation}

\z where $K,K_a,K_s$ are constants independent of $\xi$ and 

\begin{equation}\label{notatii1}
R=\frac{\rho^{1-a(q+1)}}{q+1}\ ,\ \ \lambda=iZP_0-t
\end{equation}
$$\ \ \ \ \ \ -d=\max\{ [1-a(q+1)](-2-A)\, ,\, 1-a(q+2)\}$$
\begin{equation}\label{notatii2}
\gamma={a\left(\frac{q+1}{2}-\frac{1}{2m}\right)-\frac{1}{4}}\ \ ,\ \
-\beta=[1-a(q+1)][-1-A]
\end{equation}

\end{Lemma}

The proof is given in \S\ref{prfasygminus}. 

The solution $g_-$ is a linear combination of the solutions $g_a,g_s$:

\begin{equation}\label{firstmatch}
g_-(v)=C_ag_a(v)+C_sg_s(v)
\end{equation}

\z where the coefficients $C_a$,$C_s$ depend on $\xi$.

The representations (\ref{aga}), (\ref{ags}),
(\ref{agm}) are used to find sharp estimates of $C_a$, $C_s$ by
identification of the coefficients of the independent functions whose
leading behavior is $e^{it}$, respectively $e^{-it}$:

\begin{Lemma}\label{asygas}
  
The constants $C_a,C_s$ of (\ref{firstmatch}) are analytic in $\xi$
for large $\xi$ in a strip containing $\RR_+$ and satisfy 

\begin{equation}\label{Ca}
C_a=\tilde{K}_{a}\, \xi^{\gamma+\beta}\left[
 -e^{iP_0\xi}+\frac{i}{\sqrt{2}}e^{i\pi(1+B)}e^{-iP_0\xi}+O\left(\xi^{-d+\beta}\right)\right]\ \ 
(\xi\rightarrow\infty)
\end{equation}

\begin{equation}\label{Cs}
C_s=\tilde{K}_{s}\, \xi^{\gamma+\beta}\left[
  e^{iP_0\xi}-\frac{i}{\sqrt{2}}e^{i\pi(1+A)}e^{-iP_0\xi}
+O\left(\xi^{-d+\beta}\right)\right]\ \ 
(\xi\rightarrow\infty)
\end{equation}

\z where

$$\tilde{K}_{a,s}=\frac{K}{\left(e^{i\pi A}-e^{i\pi B}\right)K_{a,s}}$$

\end{Lemma}

The proof is given in \S\ref{prfasygas}.

\subsection{Existence of a bounded solution of (\ref{eqfx}) for $m$ even}\label{matchthem}

It is enough to find values of the parameter $\xi$ for which that
$g_-$ is a multiple of $g_a$. Such values exist indeed:

\begin{Proposition}\label{Csis0}
There exist infinitely many values of $\xi$ for which $C_s=0$
\end{Proposition}

{\em{Proof}}

For large $\rho$ equation $C_s(\xi)=0$ gives, up to lower order terms, (see (\ref{Cs}))

\begin{equation}\label{valxim}
e^{i\rho P_0}-i2^{-1/2} e^{-i\rho
    P_0+i\pi(1+A)}=0
\end{equation}

\z having the solutions 

\begin{equation}\label{xiM}
\xi_M=\frac{\pi}{P_0}\left( M+\frac{1}{2m(q+1)}\right)
-\frac{i}{4P_0}\ln 2\ \ \ ,\ M\in \ZZ
\end{equation}

\z We will show that for $M$ large enough, near each $\xi_M$ there are
values of $\xi$ for which $C_s=0$.  Note that the assumption
(\ref{orderz}) is satisfied for $\xi=\xi_M$ if $M$ is large enough.

The essence of the rest of the proof is to show that an
``approximate'' solution implies the existence of an actual solution
using the principle of variation of the argument.

Let $j\in\{1;2\}$ and $\alpha>0$ with $\alpha<d-\beta$.  Consider the
circles $\mathcal{C}_j$ centered at $\xi_M$, of radius $jM^{-\alpha}$:

\begin{equation}\label{circlej}
\xi\in\mathcal{C}_j\ \ \ {\mbox{if}}\ \ \
\xi=\xi_M+jM^{-\alpha}\theta\ \ \ {\mbox{where}}\ \ \ |\theta|=1
\end{equation}

Note that for $\xi$ in the interior of the circle $\mathcal{C}_j$ we
have $\xi=O(M)$ and $M=O(\xi)$. 

Assume, to get a contradiction, that $C_s(\xi)\ne 0$ for $\xi$ in the
interior of the circle $\mathcal{C}_2$. 

Using (\ref{Cs}), (\ref{valxim}), (\ref{circlej}) for
$\xi\in\mathcal{C}_j$ we get 

\begin{equation}\label{osetimcs}
C_s(\xi)=
K_0e^{iM\pi}jM^{\gamma+\beta-\alpha}(\theta+o(1))\ \ \ (M\rightarrow +\infty)
\end{equation}

\z where $K_0$ does not depend on $M$.

So $|1/C_s(\xi) |\leq 1/(2|K_0|M^{\gamma+\beta-\alpha})(1+o(1))$ for $\xi$ on
$\mathcal{C}_2$, hence also for $\xi$ in the interior of this circle.
Then $|C_s(\xi) |\geq 2|K_0|M^{\gamma+\beta-\alpha}(1+o(1))$ in the interior of
$\mathcal{C}_2$.

However, using (\ref{osetimcs}) for $j=1$ we get that
$|C_s(\xi) |\leq |K_0|M^{\gamma+\beta-\alpha}(1+o(1))$ on the circle
$\mathcal{C}_1$, which is a contradiction.

\subsection{Matching in the case $m$ odd, $q$ even}\label{matchodd}

Denote by $\phi(v)$ the solution of (\ref{eqgv}) which equals $g_-$ on
$d$. For $v\in d$ we have $\phi(v)=C_ag_a(v)+C_s g_s(v)$ with
$C_a,C_s$ given by (\ref{Ca}), (\ref{Cs}). After analytic continuation
of $\phi(v)$ on a small path rotating around $v=0$ by an angle $m\pi$ the
solution $\phi$ becomes $\phi(e^{im\pi}v)=C_ag_a(e^{im\pi}v)+C_s
g_s(e^{im\pi}v)$ which decomposes along the small solution, and the
large solution on $s(-d)$: $\phi(e^{im\pi}v)\equiv
D_+g_{+;-\infty}(-v)+D_-g_{-;-\infty}(-v)$ for some constants
$D_+,D_-$; $g_{\pm;-\infty}$ are given in Lemma \ref{L1min}
(ii). Hence

\begin{equation}\label{decgasgpmmin}
 C_ag_a\left( e^{im\pi}v\right)+C_s g_s\left( e^{im\pi}v\right)\equiv
 D_+g_{+;-\infty}(-v)+D_-g_{-;-\infty}(-v)
\end{equation}

Note that $g_a\left( e^{im\pi}v\right)=g_a(-v)=g_a(v)$ since $g_a$ is
analytic and even, and that $g_s\left(
  e^{im\pi}v\right)=-v^{1/m}g_{a;s}(v)$ since $g_{a;s}$ is
analytic and even.  Relation (\ref{decgasgpmmin}) becomes 

\begin{equation}\label{decmin}
 C_ag_a\left( v\right)-C_s g_s\left( v\right)\equiv
 D_+g_{+;-\infty}(-v)+D_-g_{-;-\infty}(-v) 
\end{equation}

\z for $v\in d$, hence for all $v$.

We first need asymptotic estimates for the solutions
$g_a,g_s,g_{-;-\infty}$ for $v=-v_M^-$:

\begin{Lemma}\label{asygminmin}
  
  For $v=-v_M^-$ where $v_M^-$ is defined by (\ref{smmin}),
  (\ref{vmmin}) we have

\begin{equation}\label{agamin}
g_{a}(v)=K_{a}\xi^{-\beta}\left[e^{i(R-\lambda)}+e^{i\pi(1+A)}e^{-i(R-\lambda)}\right]+O(\xi^{-d})
\end{equation} 

\begin{equation}\label{agsmin}
g_{s}(v)=K_{s}\xi^{-\beta}\left[e^{i(R-\lambda)}+e^{i\pi(1+B)}e^{-i(R-\lambda)}\right]+O(\xi^{-d})
\end{equation} 

\z and
 
\begin{equation}\label{agmmin}
g_{\pm;-\infty}(v)=K^-\xi^\gamma \left[ e^{\pm i\pi/8}
e^{\pm iP_0\xi}e^{\pm i(R-\lambda)}+O(\xi^{-1+a(q+1)})\right]
\end{equation}

\z where the notations are those of (\ref{notatii1}),(\ref{notatii2}).

\end{Lemma}

The proof is given in \S\ref{prfasgminmin}. 

Denote
$$-c=-d+\beta\ \ ,\ \ -c_1=-1+a(q+1)$$

For $v$ satisfying (\ref{vmmin}) the RHS of (\ref{decmin}) becomes, in
view of (\ref{agamin}), (\ref{agsmin}), (\ref{Ca}), (\ref{Cs})

\begin{multline}\nonumber
C_ag_a(v)-C_sg_s(v)=\\\nonumber
K\xi^\gamma e^{i(R-\lambda)}\left[ -2e^{i\xi
    P_0}+\frac{i}{\sqrt{2}}\left( e^{i\pi(1+A)}+e^{i\pi(1+B)}\right)e^{-i\xi
    P_0}+O\left(\xi^{-c}\right) \right] \\\nonumber
+  K\xi^\gamma e^{-i(R-\lambda)}\left[ -e^{i\xi
    P_0}\left( e^{i\pi(1+A)}+e^{i\pi(1+B)}\right)+{i}{\sqrt{2}}
  e^{i\pi(1+A)} e^{i\pi(1+B)}e^{-i\xi P_0}\right.\\\nonumber
\left. +O\left(\xi^{-c}\right)\right]  \\\nonumber
+ \frac{K}{K^-} g_{+;-\infty}\left[ -2+\frac{i}{\sqrt{2}}\left(
       e^{i\pi(1+A)}+e^{i\pi(1+B)}\right)e^{-2i\xi
       P_0}\right. \\\nonumber\left. +\,
   O\left(\xi^{-c}\right)\right] \,
 \left(1+O\left(\xi^{-c_1}\right)\right)
\\\nonumber
   +\frac{K}{K^-}g_{-;-\infty}\left[ -e^{2i\xi P_0}\left(
       e^{i\pi(1+A)}+e^{i\pi(1+B)}\right)+{i}{\sqrt{2}} e^{i\pi(1+A)}
     e^{i\pi(1+B)}+O\left(\xi^{-c}\right)\right]\\ \nonumber\times
   \left(1+O\left(\xi^{-c_1}\right)\right)
\end{multline}

\z which in view of (\ref{decmin}) gives

\begin{multline}\label{Dp}
  D_+=K/{K^-}\left[
    -2+\frac{i}{\sqrt{2}}\left(e^{i\pi(1+A)}+e^{i\pi(1+B)}\right)
    e^{-2i\xi P_0}+O\left( \xi^{-c_1}\right)\right. \\ \left. +O\left(\xi^{-c}\right)
  \right]\, \left(1+O\left(\xi^{-c_1}\right)\right)
\end{multline}

The existence of a bounded solution of (\ref{eqgv}) on $v(-\RR)$,
hence of (\ref{eqfx}) on $\RR$ is assured if, for some value of the
parameter $\xi$ we have $D_+=0$. Indeed, there are such values:

\begin{Proposition}\label{L10p9}

There are infinitely many values of $\xi$ for which $D_+=0$,
hence for which equation (\ref{eqfx}) has a bounded solution on $\RR$.

\end{Proposition}

{\em{Proof}}

We need to find valued for $\xi$ satisfying (\ref{orderz}) such that
$D_+=0$. Up to smaller terms, $\xi$ satisfies

$$-2+\frac{i}{\sqrt{2}}\left(
       e^{i\pi(1+A)}+e^{i\pi(1+B)}\right)e^{-2i\xi
       P_0}=0$$

\z which holds for

\begin{equation}\label{xiMmodd}
\xi_M=-\frac{i}{2P_0}\ln\left(2^{-1/2}\cos\left(\frac{\pi}{2m(q+1)}\right)\right)
+\frac{\pi}{2P_0}(2M+1)
\end{equation}

In a neighborhood of each $\xi_M$ there is a zero of the function
$D_+$ if $M$ is large enough, the details being as in the proof of
Proposition \ref{Csis0}.

\section{Proof of auxiliary Lemmas}

\subsection{{Some estimates}}\label{Someestim}

{\em{I.}} For $v$ such that $s(v)$ is on $\ell$, or on $\ell^-$ we have
$|P(v)|>c_0>0$.

Indeed, note first that if $Z=0$ (so $\xi=\rho\in\RR_+$) the part $\ell_1$ of
$\ell$ corresponds in the $v$-plane to $[1,+\infty)$, the part
$\ell_{12}$ corresponds to a small half-circle avoiding $v=1$, and the
rest of $\ell$ corresponds to a subsegment of $(0,1)$. Since $P(v)\ne
0$ for $v\in\RR_+$, statement {\em{I}} follows.

In the general case when $Z\ne 0$ the picture is similar. A
continuous deformation of $s(d)$ to $\ell$ through the right
half-plane corresponds in the $v$-plane to a counterclockwise rotation
of $d$, of angle close to $\pi/(2(q+2))$. For large $v$ we have $s\sim
i(\rho+iZ)v^{q+2}/(q+2)$ so if initially $v\in d$, after the
deformation $\arg v\sim -(q+2)^{-1}\arg(\rho+iZ)$ so $\arg v$ is
small. For $v=O(1)$ ($\rho\rightarrow +\infty$), $Z$ is negligible in
(\ref{sdev}), (\ref{orderz}) so $v$ is close to $\RR_+$ and $s=0$
corresponds to $v=1$. At this point the branch of $v(s)$ (on $\ell$)
changes; denote this branch by $v_+(s)$.

{\em{II.}} We have $|v_+(s)|<1$. Indeed, for $s=i\rho\lambda P_0$
with $\lambda\in[0,1]$ substituting $v=1+\Delta$ in (\ref{sdev}) we
get $\Delta\sim \pm\sqrt{2\lambda P_0}(1-iZ\rho^{-1}/2)$ (for large
$\rho$). Then $v_+(s)\sim 1-\sqrt{2\lambda P_0}(1-iZ\rho^{-1}/2)$ and
the estimate follows.

{\em{Remark.}} From {\em{I}} and {\em{II}} it follows that $s=s_0$
is not a ramification point of $v(s)$ on the Riemann sheet of
$\CC\setminus\{0;s_0\}$ where the homotopic deformation of $s(d)$ to
$\ell$ was done; hence the presence of the path $\ell_{12}$ is not
necessary and can be omitted in the proofs.

{\em{III.}} Using {\em{I}} it follows that 
$|{vP'(v)}/{P(v)}|<const\, $, $|{v^2P''(v)}/{P(v)}|<const\, $ so that
$|\tilde{R}(v)|<const$.

\subsection{Proof of Lemma \ref{L1}}\label{ProofL1}

(i) Consider the operator $\mathcal{J}$ of (\ref{defJ}) defined on the
Banach space of functions $h$ which are continuous and bounded on
$\ell$, in the sup norm, and analytic in $\xi$.  Write
$\mathcal{J}=\mathcal{J}_+-\mathcal{J}_-$ where

$$\mathcal{J}_\pm(h)(s)= B_\pm (s)\int_{+i\infty}^s B_\mp (\sigma)\frac{\sigma^{1/2}}{W_0}R(\sigma)h(\sigma)\,
d\sigma$$

\z Using (\ref{Rdef}) and (\ref{sdev}) we have 

$$\mathcal{J}_\pm(h)(s)=const\, B_\pm (s)\int_{+i\infty}^s B_\mp
(\sigma)\frac{\sigma^{-1/2}P(v(\sigma))}{\xi
  v(\sigma)^{2q+2}}\tilde{R}(\sigma)h(\sigma)\, d\sigma$$

Note that (see \S\ref{Someestim})

\begin{equation}\label{estimR}
 \left|\tilde{R}(s)\right|\leq const\ \ \ {\mbox{for}}\ s\in\ell
\end{equation}

Let
$\|h\|_{\ell_k}=\sup_{s\in\ell_k}|h(s)|$.
 Since $B_\pm(s)\sim s^{-1/4}e^{\pm s}$ for $s\in i\RR_+$ and the
solutions of (\ref{eqB}) are bounded at $s=0$ then $|B_\pm(s)|<const$
for $s\in\ell_2$ and $|B_\pm(s)|<const\, |s|^{-1/4}$ for $s$ in
$\ell_1$, $\ell_{12}$, $\ell_3$ or $\ell_4$.

For $s\in\ell_1$ we have $|v(s)|\geq 1/2$.
Also, since $|(v-1)^2P(v)/|v|^{q+2}|<const$ on $v(\ell_1)$ then
$|v|^{-q-2}<const\, \rho |s|^{-1}$. 

Hence (using also $|P(v)/v^q|<const$)

$$\left| \int_{i\rho P_0/2}^{+i\infty}B_\pm(\sigma)\frac{\sigma^{-1/2}P(v(\sigma))}{\xi
  v(\sigma)^{2q+2}}\tilde{R}(\sigma)h(\sigma)\, d\sigma\right|<const\,\|h\|_{\ell_1}
\int_{\rho P_0/2}^{+\infty}\rho^{-1}\tau^{-3/4}\, d\tau$$
$$=const\, \rho^{-3/4}\|h\|_{\ell_1}$$

For $s\in\ell_2$ we have $const_1>|v|>const_2 >0$ hence

$$\left| \int_{i\rho P_0/2}^{i\rho P_0/2}B_\pm(\sigma)\frac{\sigma^{-1/2}P(v(\sigma))}{\xi
  v(\sigma)^{2q+2}}\tilde{R}(\sigma)h(\sigma)\, d\sigma\right|<const\, \|h\|_{\ell_2}
\int_0^{\rho P_0/2} \frac{\tau^{-3/4}}{\rho}\, d\tau$$
$$=const\, \rho^{-3/4}\|h\|_{\ell_2}$$

If $\ell_{12}$ is part of $\ell$ the possible integral on the
segment $[i\rho P_0, i\xi P_0]$ does not modify these estimates since
the length of this segment is $O(1)$ as $\rho\rightarrow+\infty$.

For $s\in\ell_3$ the
corresponding $v$ satisfies $const\,\rho^{-a}<|v|<const\, <1$ (since
the branch of $v(s)$ is $v_+(s)$, see \S\ref{Someestim} {\em{II}}). Then
substituting $\sigma=i\rho\tau$ in

$$T:=\left| B_\pm(s)\int_{i\rho P_0/2}^{s_M}B_\mp
  (\sigma)\frac{\sigma^{-1/2}P(v(\sigma))}{\xi
    v(\sigma)^{2q+2}}\tilde{R}(\sigma)h(\sigma)\, d\sigma\right|$$

\z and (since $|P(v)|<const$) we get

\begin{equation}\label{oformula}
T <const\, \rho^{-1} \|h\|_{\ell_3}\, \int_{1/4}^{1/2(1-\rho^{-a})}
\frac{\tau^{-3/4}}{|v(i\rho\tau)|^{2q+2}}\, d\tau
\end{equation}

\z where now $v$ is given implicitly by $\tau=(1+z)(v-1)^2P(v)$, with the
notation $z=iZ\rho^{-1}$.

Let $\tilde{v}(s)$ denote $v(s)$ for $z=0$.
We have $\tilde{v}_+\in [c_0\rho^{-a},v_0]\subset (0,1)$. Then
$|v_+(s)-\tilde{v}_+(s)|\leq |z| \sup_{|z|<const\,\rho^{-1}}|v'_+(z)|< const\,
|z|\rho^{-aq}$. It follows that $|v|^{-1}<const\, \tilde{v}_+^{-1}$ on
$\ell_3$.

Then (\ref{oformula}) is less than

$$const\, \rho^{-1} \|h\|_{\ell_3}\, \int_{1/4}^{1/2(1-\rho^{-a})}
\frac{1}{\tilde{v}(\tau)^{2q+2}}\, d\tau = const\, \rho^{-1} \|h\|_{\ell_3}\,
\int_{K_0\rho^{-a}}^{v_0}\frac{1-\tilde{v}}{v_1^{q+2}}\, d\tilde{v} $$
$$\leq const\, \rho^{-1+a(q+1)}\| h\|_{\ell_3}$$

Finally, for $s\in \ell_4$ we have $s=O(\rho)$,
$|s_M-s(v_{P})|<const\, \rho^{1-a(q+1)}$ and $|v|>const\, \rho^{-a}$
 so that 

 $$\left| B_\pm(s)\int_{s_M}^s B_\mp
   (\sigma)\frac{\sigma^{-1/2}P(v(\sigma))}{\xi
     v(\sigma)^{2q+2}}\tilde{R}(\sigma)h(\sigma)\,
   d\sigma\right|<const\, \rho^{1-a(q+1)}\|h\|_{\ell_4}$$

Therefore 

$$\left|\mathcal{J}(h)(s)\right|<const\,\max\{\rho^{-3/4},\rho^{-1+a(q+1)}\}
\|h\|$$

\z and $\mathcal{J}$ is a contraction for $a$ satisfying
(\ref{orderz}). Equation (\ref{eqhmin}) has a unique solution
$h_-=(1-\mathcal{J})^{-1}\mathcal{J}(B_-)$.

To show (\ref{ttxxvp}) we write
$h_-=B_-+\left( 1-\mathcal{J}\right)^{-1}\mathcal{J}(B_-)$. We showed that

$$|\mathcal{J}B_-(s)|\leq const\, \rho^{-3/4}
\|B_-\|_{\ell_1\cup\ell_{12}\cup\ell_2}\, +\, 
  const\, \rho^{-1+a(q+1)}\| B_-\|_{\ell_3}$$

\z and since $\| B_-\|_{\ell_3}< const\, \rho^{-1/4}$ and
$\|(1-\mathcal{J})^{-1}\|<const$, the estimate (\ref{ttxxvp}) follows.

(ii) is immediate since $B_-$, therefore $h_-$ are exponentially small
on $s(d)$.

Finally, (iii) holds since $B_\pm(s)$ do not depend on $\xi$ and the operator
$\mathcal{J}$ is analytic in $\xi$.

\subsection{Proof of Lemma \ref{L1min}}\label{ProofL1m}

Denote by $\ell^-_1$ the segment of $\ell^-$ from $+i\infty$ to
$2i\rho P_0$, by $\ell^-_3$ the segment from $2i\rho P_0$ to
$s_M^-$, and by $\ell^-_4$ the segment from $s_M^-$ to $s_F^-$

The estimates on $\ell^-_1$ (respectively $\ell^-_3$) are the same
as in the proof of Lemma \ref{L1}(i) on $\ell_1$ (respectively
$\ell_3$). For the integral on $\ell^-_4$ the integrand has the same
bounds as on $\ell^-_3$, but since the path of integration has an
$O(1)$ length the contribution of this term does not modify the estimate.

\subsection{Proof of Lemma \ref{L2}}\label{ProofL2}

\subsubsection{The integral equation}\label{Wnudepxi}

The substitution $y=i\xi
v^{q+1}/(q+1)$ transforms (\ref{eqgvsmall}) to

\begin{equation}\label{eqhy}
h''(y)+\frac{r}{y}h'(y)-h(y)=\left(v(y)^2-2v(y)\right)h(y)
\end{equation}

\z where $r$ is given by (\ref{defr}).

The LHS of (\ref{eqhy}) is a Bessel equation, with an analytic solution

\begin{equation}
\tilde{h}_a(y)=\int_{-1}^1e^{-py}(1-p^2)^A\, dp
\end{equation}

\z and the branched solution 

\begin{equation}
\tilde{h}_s(y)=y^{1-r}{\tilde{h}}_{a;s}(y)\ \ ,\ \
{\mbox{where}}\ {\tilde{h}}_{a;s}(y)= \int_{-1}^1e^{-py}(1-p^2)^B\, dp
\end{equation}

\z Their Wronskian is $W[\tilde{h}_a,\tilde{h}_s]=const\, y^r$.

Going back to the variable $v$, let
$\tilde{g}_{a,s}(v)=\tilde{h}_{a,s}(y(v))$. Then

\begin{equation}\label{Wgags}
W[\tilde{g}_a,\tilde{g}_s]=const\, \xi^{1-r} v^{-1+1/m}
\end{equation}

The integral form of equation (\ref{eqgvsmall}) is then (\ref{ginteq})
where

$$\mathcal{G}(g)(v)=\int_0^v \left[
  \tilde{g}_a(v)\tilde{g}_s(t)-\tilde{g}_s(v)\tilde{g}_a(t)\right]W[\tilde{g}_a,\tilde{g}_s]^{-1}
\xi^2 t^{2q+1}(2-t)g(t)\,dt$$

\z which using (\ref{Wgags}) yields (\ref{defcalG}).

\subsubsection{Estimates}

We show that the operator $\mathcal{G}$ is contractive in the Banach space of
continuous and bounded functions on $\mathcal{S}$ (with the $\sup$ norm).

We need the following estimates:

Consider first $v\in \mathcal{S}_1$. We have $y\equiv y(v)\in
i[0,\rho^{1-a(q+1)}/(q+1)]\subset i\RR_+$. Note that 
$\tilde{h}_a,\tilde{h}_s$ are bounded for small $y$; also, for large $y$
on the segment $y(v)$ we have $|\tilde{h}_{a,s}(y)|<const\, |y|^{-1-A}$
(from (\ref{estimgas}) used for $\alpha=-A$ and $\alpha=-B$). 

Denote $v_1=\kappa_0 (q+1)^{1/(q+1)}\rho^{-1/(q+1)}$. Then
$y(v_1)=i$. 

For $v\in [0,v_1]$ we have $|\tilde{g}_{a,s}(v)|\leq
\sup_{y\in [0,i]}|\tilde{h}_{a,s}(y)|=const$ so that 

$$| \mathcal{G}(g)(v)|< const\, \rho^{1+r}\sup_{v\in [0,v_1]}|g(v)|
\int_0^{\rho^{-1/(q+1)}}t^{2q+2-\frac{1}{m}}\, dt$$
$$< const\, \rho^{-1/(q+1)}\sup_{v\in [0,v_1]}|g(v)|$$

For $v\in [v_1,v_F]$
we have $|\tilde{g}_{a,s}(v)|=\tilde{h}_{a,s}(y(v))|<const\, |\rho
v^{q+1}|^{-1-A}$ so that 

\begin{multline}
  | \mathcal{G}(g)(v)|< const\, \rho^{-1/(q+1)}\sup_{v\in
    [0,v_1]}|g(v)|\\
  + const\, \rho^{-2A-\frac{1}{m(q+1)}}\sup_{v\in [v_1,v_F]}|g(v)|\,
  |v_F|^{(q+1)(-1-A)}\, \int_{\rho^{-1/(q+1)}}^{\rho^{-a}}
  t^{(1-A)(q+1)-\frac{1}{m}}\, dt\\ < const\,
  \rho^{1-a(q+2)}\sup_{v\in\mathcal{S}_1}|g(v)|
\end{multline}

\z Therefore 

\begin{equation}\label{Gcontraction}
\sup_{v\in\mathcal{S}_1 }|\mathcal{G}(g)(v)|\le const\, \sup_{v\in\mathcal{S}_1 }|g(v)|
\end{equation}

For $v\in \mathcal{S}_2$ we have  

$$\xi v^{q+1}= \xi v_F^{q+1}\left[1+
  O\left(\rho^{-a(q+1)}\right)\right]\ \ (\rho\rightarrow\infty)$$

\z hence this path yields a negligible contribution  to
(\ref{Gcontraction}) and $\left\|\mathcal{G}(g)\right\|\le const\,
\rho^{1-a(q+2)}\|g\|$.

Then $\mathcal{G}$ is contractive for $\rho$ large enough.
Taking $C_1=1,C_2=0$ in (\ref{ginteq}) (respectively $C_1=0,C_2=1$) we
get the analytic (respectively branched) solution of
(\ref{eqgv}) 

\begin{equation}\label{estega}
g_{a,s}=(I-\mathcal{G})^{-1}(\tilde{g}_{a,s})=\tilde{g}_{a,s}+
(I-\mathcal{G})^{-1}\mathcal{G} \tilde{g}_{a,s}
\end{equation}

\z and the estimate (\ref{estimgasD}) follows.

\subsection{Proof of Lemma \ref{asygminus}}\label{prfasygminus}

We first estimate the leading order $B_-$ of $h_-$.

Formula (\ref{formBpm}) defines $B_-(s)$ for $s$ between $+i\infty$
and $0$. To get the values of $B_-(s)$ in the matching region, where 

$$s(v)=i\left\{ \rho P_0
-t^{q+1}R \right\} +o\left(\rho^{1-a(q+2)}\right)
+O\left( \rho^{1-a(2q+3)}\right)\ {\mbox{with}}\ R=\frac{\rho^{1-a(q+1)}}{q+1}$$

\z we need to find the analytic continuation of (\ref{formBpm}) along
$\ell$. 

We remark that in fact $B_-$ also satisfies 

\begin{equation}\label{deformBpm}
B_-(e^{i\theta}s)=\int_{e^{-i\gamma}[1,+\infty)}e^{- ps}\left( p^2-1\right)^{-3/4}\, dp\ \
\ \ \ \mbox{for}\ \  s\in i\RR_+\ \ {\mbox{if}}\
\sin(\theta-\gamma)\leq 0
\end{equation}
 \z and formula (\ref{deformBpm}) can be used to determine the analytic
continuation of $B_-$ along $\ell$: when $s$ goes on a small circle,
clockwise around the origin, the path of integration in
(\ref{formBpm}) must be rotated by an angle of $2\pi$
(around the singular points $p=\pm 1$). Then for $s\in\ell_2$

\begin{multline}\label{ACBM}
B_-(s)=\int_1^{-1} e^{- ps} e^{-3\pi i/4}\left(1-
p^2\right)^{-3/4}\, dp\, \\ +\, \int_{-1}^1 e^{- ps} e^{-9\pi i/4}\left(1-
p^2\right)^{-3/4}\, dp
 +\, \int_1^{+\infty}e^{- ps}e^{-3\pi i}\left( p^2-1\right)^{-3/4}\, dp
\end{multline}
\begin{equation}
=2^{1/2}\int_{-1}^1e^{- ps}\left(1-
p^2\right)^{-3/4}\, dp\, - \, \int_1^{+\infty}e^{- ps}\left( p^2-1\right)^{-3/4}\, dp
\end{equation}

Using now (\ref{from6.34}) and  (\ref{estimgas}) we get

\begin{multline}
B_-(s)=const\, \rho^{-1/4} \left( ie^{-i\rho
    P_0}e^{it^{q+1}R}+2^{1/2}ie^{i\rho P_0}e^{-it^{q+1}R} \right)\\
\times\left[ 1+o\left(\rho^{1-a(q+2)}\right)
+O\left( \rho^{-a(q+1)}\right) \right]
\end{multline}

Finally, in view of Lemma \ref{L1}(ii), and noting that in the
matching region (\ref{vm}) we have $P(v)=P_0+O(\rho^{-a})$ and
$v^{-(q+1)/2+1/(2m)}=\rho^{a[(q+1)/2-1/(2m)]}(1+o(1))$ Lemma \ref{asygminus} follows.

Analyticity in $\xi$ follows as in \S\ref{ProofL1}.

\subsection{Proof of Lemma \ref{asygas}}\label{prfasygas}

{\em{I. Analytic dependence on $\xi$ of $C_a,C_s$.}}

The solutions $g_\pm(v),g_{a,s}(v)$, depend analytically on $\xi$ (Lemmas
\ref{L1}, \ref{L2}). We have $g_-(v(s))\sim const\, s^{-1/4}e^{-s}$
for $s$ in the initial part of $\ell$. After analytic continuation
along $\ell$, rotating $2\pi$ about $s=0$, the solution $g_-$ becomes
$AC g_-(v)=A(\xi)g_-(v)+B(\xi)g_+(v)$. Evaluating this relation for
two values of $v$ it follows that $A(\xi)$, $B(\xi)$ are analytic. The
same argument shows that the decomposition of $g_a,g_s$ along
$g_+,g_-$ has analytic coefficients.

{\em{II. Estimates.}}

Introducing (\ref{aga}),(\ref{ags}),(\ref{agm}) in (\ref{firstmatch})
and identifying the coefficients of $e^{i(R+\lambda)}+O(\rho^{-c})$, respectively
$e^{-i(R+\lambda)}+O(\rho^{-c})$ we get a system of two
equations for the unknowns $C_a,C_s$, whose solutions have the
form (\ref{Ca}), (\ref{Cs}).

\subsection{Proof of Lemma \ref{asygminmin}}\label{prfasgminmin}

{\em{I.}} We first note that the result in Lemma \ref{L2} can be improved: the
integral equation (\ref{ginteq}) has a unique solution on the longer
path $[0,v_M]\cup [v_M, -v_M^-]$ and the estimates (\ref{estimgasD}) hold.

Indeed, on the segment $[v_M, -v_M^-]$ we have
$y=i\rho^{1-a(q+1)}/(q+1)+\mu+O(\rho^{-a})$ where $\mu$ is bounded,
therefore the estimates of $\tilde{h}_{a,s}$ are the same as on the
segment $[0,v_M]$ (possibly after increasing  $const$).
Since the segment $[v_M, -v_M^-]$ has an $o(1)$ length, its
contribution to the estimates in the proof of Lemma \ref{L2} is
negligible.

{\em{II.}} Note that 

$$s\left( -v_M^-\right)=i\rho P_0
-i\frac{\rho^{1-a(q+1)}}{q+1}+i\left( 2iZP_0-t\right)+o(1)\ \ \
(\rho\rightarrow +\infty)$$

The rest follows in a straightforward way from (\ref{estimgasD}),
(\ref{ttxxvm}) and Remark \ref{aint1inf}(iv).

\section{Appendix}

\subsection{The paths $s(d)$, $s(-d)$}\label{s(d)}

Denote $\theta=\frac{\pi}{2(q+2)}$.

Let $v\in d$ (respectively $v\in -d$), so
$v=e^{-i\theta}t$ with $t\in\RR_+$, (respectively
$t\in\RR_-$).  Then

$$\Re[s(v)]=\rho\left(
  \frac{1}{q+2}t^{q+2}- \frac{\cos \theta}{q+1}t^{q+1}\right)
+Z\left(\frac{\sin\theta}{q+1}t^{q+1}+P_0\right) $$

\z and

$$\Im[s(v)]=Z\left(
  \frac{1}{q+2}t^{q+2}- \frac{\cos\theta}{q+1}t^{q+1}\right)
+\rho\left(-\frac{\sin\theta}{q+1}t^{q+1}+P_0\right)$$

For large $|t|$ we have $\Re s(v)\sim \rho t^{q+2}/(q+2)$ and $\Im
s(v)\sim Z/(q+2)t^{q+2}$. Hence $s(d)$ (and $s(-d)$ if $q$ is even)
comes from infinity in right half-plane, asymptotic to a line of small
slope ($Z\rho^{-1}$). 

To obtain the form of $s(d)$ consider first the case $Z<0$. Then
$s(d)$ goes to infinity in the fourth quadrant.  Furthermore, we have
$\Im s(v)=0$ only for $t=[(q+2)\sin\theta]^{-1/(q+1)}(1+o(1))$ (and at
this point $\Re s(v)<0$); also, $\Re s(v)=0$ only for
$t=(q+2)/(q+1)(\cos\theta)(1+o(1))$ (at this point $\Im s(v)<0$) and
for $t=[-Z/(q+2)/ \cos\theta \rho^{-1}]^{1/(q+1)}(1+o(1))$.  Therefore
$s(d)$ comes from $\infty$ in the fourth quadrant, crosses $i\RR_-$,
then $\RR_-$, $i\RR_+$ and ends at $s_0$ in the first quadrant.

The form of $s(d)$ in the case $Z>0$ is obtained similarly. The
curve comes from infinity in the first quadrant, crosses $\RR_+$, then
$i\RR_-$, and $\RR_-$ and ends in the second quadrant at $s_0$.

For $q$ even, $s(-d)$ starts at $s_0$ (in the second quadrant if
$Z>0$, respectively in the first quadrant if $Z<0$) crosses $i\RR_+$
(respectively $\RR_+$), then goes to infinity in the first
(respectively fourth) quadrant.

\subsection{Asymptotic formulas}\label{form6.34}

\begin{Remark}\label{aint1inf}

Let $T>0$ and $\alpha\in (0,1)$. 

(i) We have

\begin{equation}\label{6.34}
\int_0^{+\infty}e^{-Tq}q^{-\alpha}(2\pm iq)^{-\alpha}\, dq =
 2^{-\alpha}\Gamma(1-\alpha)T^{-1+\alpha}\left(1+O\left(T^{-1}\right)\right)
\end{equation}

(ii) Also

\begin{multline}\label{from6.34}
  \Phi_\pm(T)\equiv \int_1^{+\infty}e^{\pm ipT}\left(p^2-1\right)^{-\alpha}\, dp\\
  =e^{\pm i\pi (1-\alpha)/2}2^{-\alpha}\Gamma(1-\alpha)
  e^{\pm iT}T^{-1+\alpha}\left( 1+O(T^{-1})\right)
\end{multline}

(iii) Denote

\begin{equation}\label{IdeT}
I(T)=\int_{-1}^1 e^{-ipT}(1-p^2)^{-\alpha}\, dp
\end{equation}

\z We have 
\begin{multline}\label{estimgas}
I(T)=\left[ e^{iT-i\pi (1-\alpha)/2}+e^{-iT+i\pi(1-\alpha)/2}\right]
2^{-\alpha}\Gamma(1-\alpha)T^{-1+\alpha}\left(1+O\left(T^{-1}\right)\right)
\end{multline}

(iv) Let also $K\in\CC$ independent of $T$ and
$\tilde{T}=O(T^{-\gamma})$ ($T\rightarrow +\infty$) where
$\gamma>0$. Then the integral (\ref{IdeT}) satisfies

\begin{gather}\label{estIDKT}
I\left( T+K+\tilde{T}\right)=\left[ e^{i(T+K)-i\pi
    (1-\alpha)/2}+e^{-i(T+K)+i\pi(1-\alpha)/2}\right]\\ \nonumber
\times\, 2^{-\alpha}\Gamma(1-\alpha)T^{-1+\alpha} 
\left(1+O\left(T^{-1}\right)+O\left(T^{-\gamma}\right)\right)
\end{gather}

\end{Remark}

{\em{Proof.}}

(i) Formula (\ref{6.34}) is a direct consequence of Watson's Lemma \cite{Wasow}.

(ii) To show (\ref{from6.34}) we substitute $r=p-1$ then rotate the
path of integration in the formula defining $\Phi_\pm$ by an angle
$\pm\pi/2$; it follows that

\begin{multline}
\Phi_\pm (T)=e^{\pm iT}\int_0^{-i\infty}e^{\pm irT}r^{-\alpha}(r+2)^{-\alpha}\,
dr\\= \pm i e^{\mp i\pi\alpha/2}e^{\pm iT}\int_0^{+\infty} e^{-\tau
  T}\tau^{-\alpha}(2-i\tau)^{-\alpha}\, d\tau
\end{multline}

\z and (\ref{6.34}) yields (\ref{from6.34}).

(iii) The path of integration in (\ref{estimgas}) can be deformed to two
vertical half-lines in the half-plane $\Im p<0$ (where the branches of
the power are determined by analytic continuation of the usual branch
on $(-1,1)$): 

\begin{gather}\nonumber
I(T)=\left(\int_{-1}^{-1-i\infty}+\int_{1-i\infty}^1\right)
e^{-ipT}(1-p^2)^\alpha\, dp\\\nonumber
=-ie^{i\pi\alpha/2} e^{iT}\int_0^{+\infty}e^{-sT} s^{-\alpha}(2+is)^{-\alpha}\,
  ds\ \ \ \ \ \ \ \ \ \ \ \ \ \ \ \ \ \ \  \\\nonumber\ \ \ \ \ \ \ \ \ \ \ \ \ \ \ \ \ \ \ \ 
+\, i e^{-i\pi\alpha/2}e^{-iT}\int_0^{+\infty}e^{-sT}
  s^{-\alpha}(2-is)^{-\alpha}\, ds
\end{gather}

\z and using (\ref{6.34}) the estimate (\ref{estimgas}) follows.

(iv) For large $T>0$ and $p$ with $\Re p\in [-1,1]$, $\Im p< 0$, and
$-\Im p$ large enough we have $\Re
[-ip(T+K+\tilde{T})]=[T+\Re(K+\tilde{T})]\Im p +[\Im (K+\tilde{T})]\Re
p<0$. Then the path of integration in (\ref{IdeT}) can be deformed
as in the proof of (iii), yielding

\begin{multline}
  I\left( T+K+\tilde{T}\right)=-ie^{i\pi\alpha/2}
  e^{i(T+K+\tilde{T})}\int_0^{+\infty}e^{-s((T+K+\tilde{T}))}
  s^{-\alpha}(2+is)^{-\alpha}\,
  ds\\
  +i
  e^{-i\pi\alpha/2}e^{-i(T+K+\tilde{T})}\int_0^{+\infty}e^{-s(T+K+\tilde{T})}
  s^{-\alpha}(2-is)^{-\alpha}\, ds
\end{multline}

\z and using again Watson's Lemma (\ref{estIDKT}) follows.
\vskip 1cm

\begin{center}{\bf{ACKNOWLEDGMENTS}} \end{center}

We would like to thank Fran\c cois Treves for
suggesting this problem to us and for very interesting discussions.

\end{document}